\documentclass{article}

\usepackage{amsfonts}
\usepackage{amsmath}
\usepackage{amssymb}
\usepackage[all,dvips]{xy}
\usepackage{graphics}
\usepackage{cite}

\newcommand{\Pbar}{{\ensuremath{\overline{P}}}}

\newcommand{\rad}{\ensuremath{\textup{rad}\,}}
\newcommand{\Hom}{\ensuremath{\textup{Hom}}}

\newcommand{\Ext}{\ensuremath{\textup{Ext}}}
\newcommand{\End}{\ensuremath{\textup{End}}}
\newcommand{\Ann}{\ensuremath{\textup{Ann}}}

\newcommand{\pd}{\ensuremath{\textup{pd}\,}}
\newcommand{\id}{\ensuremath{\textup{id}\,}}

\newcommand{\dcc}{\Ext^2_C(DC,C)}

\newcommand{\longto}{\longrightarrow}

\newcommand{\zG}{\ensuremath{\Gamma}}
\newcommand{\zS}{\ensuremath{\Sigma}}

\newcommand{\za}{\ensuremath{\alpha}}
\newcommand{\zb}{\ensuremath{\beta}}
\newcommand{\ze}{\ensuremath{\epsilon}}
\newcommand{\zd}{\ensuremath{\delta}}
\newcommand{\zg}{\ensuremath{\gamma}}

\newcommand{\zL}{\ensuremath{\Lambda}}
\newcommand{\zs}{\ensuremath{\sigma}}

\newtheorem{thm}{Theorem}[section]

\newtheorem{lem}[thm]{Lemma}
\newtheorem{example}[thm]{Example}
 
\newtheorem{remark}[thm]{Remark}

\newenvironment{pf}{{Proof}.}

\newcommand{\qed}{\ensuremath{\,\hfill\Box}}

\begin{document}
\title{On one-point extensions of cluster-tilted algebras}
\author{Miki Oryu and Ralf Schiffler}
\maketitle
\begin{abstract}
We define an operation which associates to a pair $(B,M)$ where $B$ is a cluster-tilted algebra and $M$ is a  $B$-module which lies in a local slice of $B$, a new cluster-tilted algebra $B'$. In terms of the quivers, this operation corresponds to adding one vertex (and arrows). 
\end{abstract}
 \section{Introduction} 
 Cluster-tilted algebras are finite dimensional associative algebras which were introduced in the context of categorifications of cluster algebras in \cite{CCS,BMR}. Since then, these algebras have been the subject of many studies, see for example 
\cite{ABS,ABS2,ABS3,ABS4,AsRe,BFPPT,BOW,BBT,BMR2,CCS2,KR}.
  
 A striking property shown in \cite[Theorem 2.13]{BMR4} is that deleting a vertex of a cluster-tilted algebra produces again a cluster-tilted algebra, more precisely:

\emph{If $B$ be a cluster-tilted algebra with quiver $Q_B$ and   $e_x$ is a primitive idempotent corresponding to a vertex $x$ of $Q_B$, then $B/Be_xB$ is cluster-tilted.}

Note that the quiver of $B/Be_x B$ is obtained from $Q_B$ be deleting the vertex $x$ and all arrows incident to it.

The goal of this paper is to give a construction which is inverse to the deletion of vertex. Given a cluster-tilted algebra $B$ with quiver $Q_B$ we want to construct a cluster-tilted algebra $B'$ whose quiver contains $Q_B$ as a full subquiver and has one additional vertex $x$, the \emph{extension vertex}, such that  $B=B'/B'e_xB'$.

Our first construction  is the one-point extension $B[P]$ of the algebra $B$ with respect to a projective B-module $P$. We prove that $B[P]$ is cluster-tilted if there exists a local slice $\zS$ in $\textup{mod}\,B$ that contains $P$, see Theorem \ref{thm main}. In this situation, the extension vertex is a source in the quiver of $B[P]$. We also give an example which shows that the condition that $P$ lies on a local slice is necessary, see Example \ref{exnonex}.

Dually, the one-point coextensions of $B$ with respect to an injective module is cluster-tilted if the injective module lies in a local slice. 
In this situation, the extension vertex is  a sink.
 
 Our second construction yields cluster-tilted algebras in which the extension vertex can have both incoming and outgoing arrows. Given any $B$-module $M$ lying on a local slice $\zS$ in $\textup{mod}\,B$, let $C=B/\Ann \zS$ be the quotient of $B$ by the annihilator of the local slice. It is known that $C$ is a tilted algebra, see \cite{ABS2}. Let $C[M]$ be its one-point extension with respect to $M$, and let $B'$ be the relation extension of $C[M]$, that is
 \[B'=C[M]\ltimes\Ext_{C[M]}(DC[M],C[M]).\]
 Then we show that $B' $ cluster-tilted such that $B'/B' e_x B' =B$, see Theorem \ref{thm main2}.
 
 Along the way, we prove that for arbitrary algebras of global dimension 2, the operation of one-point extension with respect to a projective module and  the operation of relation extension commute, see Theorem \ref{thm 1}.
 
 \section{Preliminaries} \label{sect2}
 Let $k$ be an algebraically closed field. The algebras in this paper  are always finite dimensional, basic, associative $k$-algebras and the modules are always finitely generated. If $\zL$ is a $k$-algebra, we denote by $\textup{mod}\,\zL$ the category of finitely generated right $\zL$-modules, by $\zG(\textup{mod}\,\zL)$ its Auslander-Reiten quiver and  by $\tau_\zL$ its Auslander-Reiten translation. Furthermore, $Q_\zL$ will denote the ordinary quiver of $\zL$, and $P_\zL(i), I_\zL(i)$ and $S_\zL(i)$ the indecomposable projective, injective and simple $\zL$-module at the vertex $i$ of $Q_\zL$, respectively.
Throughout the article, we use the notation $aM$ for the direct sum of $a$ copies of the module $M$.  
 For the representation theory of $k$-algebras, we refer to \cite{ARS,ASS}.
 
 \subsection{Slices and local slices}
 A \emph{path} in $\textup{mod}\,\zL$ with source $X$ and target $Y$ is a sequence of non-zero morphisms $X=X_0\to X_1\to \cdots \to X_s=Y$ where $X_i \in \textup{mod}\,\zL$ for all $i$, and $s\ge 1$.   A \emph{path} in $\zG(\textup{mod}\,\zL)$ with source $X$ and target $Y$ is a sequence of arrows $X=X_0\to X_1\to \cdots \to X_s=Y$ in the Auslander-Reiten quiver. A \emph{sectional path}
 is a path  $X=X_0\to X_1\to \cdots \to X_s=Y$ in $\zG(\textup{mod}\,\zL)$ such that for each $i$ with $0<i<s$, we
have $\tau_\zL X_{i+1}\ne X_{i-1}$.
 
 A \emph{slice} $\zS$ in $\zG(\textup{mod}\,\zL)$ is a set of indecomposable $\zL$-modules such that 
 \begin{itemize}
 \item[(S1)] $\zS$ is sincere.
 \item[(S2)] Any path in $\textup{mod}\,\zL$ with source and target in $\zS$ consists entirely of modules in $\zS$.
 \item[(S3)] If $M$ is an indecomposable non-projective $\zL$-module then at most one of $M, \tau_\zL M$ belong to $\zS$.
 \item[(S4)] If $M\to S$ is an irreducible morphism with $M$ and $S$ indecomposable and $S\in \zS$, then either $M$ belongs to $\zS$ or $M$ is non-injective and $\tau_\zL^{-1}M$ belongs to $\zS$.
 \end{itemize} 
 
 A \emph{local slice} $\zS$  in $\zG(\textup{mod}\,\zL)$ is a set of indecomposable $\zL$-modules inducing a connected full subquiver of $\zG(\textup{mod}\,\zL$) such that
 \begin{itemize}
 \item[(LS1)] If $X\in \zS$ and $X\to Y$ is an arrow in $\zG(\textup{mod}\,\zL)$ then either  $Y\in \zS$ or $\tau Y\in \zS$.
 \item[(LS2)] If $Y\in \zS$ and $X\to Y$ is an arrow  in $\zG(\textup{mod}\,\zL)$ then either  $X\in \zS$ or $\tau^{-1}X \in \zS$.
 \item[(LS3)] For every sectional path $X=X_0\to X_1\to \cdots \to X_s=Y$ in $\zG(\textup{mod}\,\zL)$ with $X$ and $Y$ in $\zS$ we have $X_i\in \zS$, for $i=0,1,\ldots,s$.
 \item[(LS4)] The number of indecomposables in $\zS$  equals the number of isoclasses of simple $\zL$-modules. 
 \end{itemize} 

\subsection{Sections and left sections} 
 Let $\zL$ be a $k$-algebra and $\zG$ be a connected component of the Auslander-Reiten quiver $\zG(\textup{mod}\,\zL)$.
A full connected subquiver $\zS$ of $\zG$  is called a \emph{section} if 
\begin{itemize}
\item[(s1)]  $\zS$ has no oriented cycles.
\item[(s2)]  $\zS$ intersects every $\tau_{\zL}$-orbit in $\zG$ exactly once.
\item[(s3)]  Each path in $\zG$ with source and target in $\zS$ lies entirely in $\zS$.
\end{itemize}

A full connected subquiver $\zS$ of $\zG$  is called a \emph{left section} if 
it satisfies condition (s1) and (s3) above as well as  the following condition.
\begin{itemize}
\item[(s2')] For any indecomposable $X$ in  $\zG$ such that there exists an indecomposable $Y$ in $\zS$ and a path in $\zG$ from $X$ to $Y$, there exists a unique $m\in \mathbb{Z}_{\ge 0}$ such that $\tau_{\zL}^{-m} X\in \zS$.
\end{itemize}

\subsection{Tilted algebras}
Let $A$ be a hereditary $k$-algebra and let $n$ be the number of isoclasses of simple  modules in  $\textup{mod}\,A$. An $A$-module $T$ is called a \emph{tilting module} if $\Ext^1_A(T,T)=0$, and $T$ is the direct sum of $n$ non-isomorphic indecomposable $A$-modules. The corresponding endomorphism algebra $\End_A T$ is called a \emph{tilted algebra}. 

Tilted algebras can be characterized by the  existence of a slice as follows.
\begin{thm}
\cite[4.2.3]{R}  Let $C=\End_{A} T $ be a tilted algebra. Then the class of $C$-modules $\Hom_{A}(T,I)$ such that $I$ is an indecomposable injective $A$-module forms a slice in $\textup{mod}\,C$. Conversely, any slice in any module category is obtained in this way.
\end{thm}

If $C$ is a tilted algebra and $M$ is a $C$-module lying on a slice in $\zG(\textup{mod}\,C)$ then the projective dimension and the injective dimension of $M$ are at most 1. Tilted algebras have global dimension at most 2. 
 
%
There is another characterization of tilted algebras in terms of  sections due to  Liu \cite{L} and Skowronsky \cite{S}. 
We will need the following result which has been shown by Assem.
\begin{thm}
 \cite[Theorem A]{A}\label{thm assem} Let $C$ be a $k$-algebra and let $\zS$ be a left section in a component $\zG$ of $\zG(\textup{mod}\,C)$ such that 
\[\Hom_{C}(\tau_{C}^{-1}E',E'' ) = 0 ,\] for all $E',E'' \in \zS$.
Then $C/\Ann_C\zS$ is a tilted algebra having $\zS$ as a slice.
\end{thm}

For further information on tilted algebras we refer the reader to \cite{ASS}.

\subsection{Relation extensions and cluster-tilted algebras}
Following \cite{ABS}, we make the following definitions.
If $C$ is an algebra of global dimension at most 2 then its relation extension $R(C)$ is defined as
\[R(C)=C\ltimes \Ext^2_C(DC,C),\]
where $D=\Hom(-,k)$ denotes the standard duality. 
The quiver $Q_{R(C)}$ has the same vertices as the quiver $Q_C$, and the arrows in  $Q_{R(C)}$ are the arrows in $Q_C$ plus one new arrow  $i\to j$ for each relation from $j$ to $i$ in a minimal system of relations for $C$, see \cite{ABS}. Moreover, the dimension of $\Ext^2_C(I_C(j),P_C(i))$ is equal to the number of nonzero paths from $i $ to $j$ in the quiver $Q_{R(C)}$ that use a new arrow. The paths in $Q_{R(C)}$ corresponding to $\dcc$ are called the \emph{new paths} in the relation extension $R(C)$. 

An algebra is called a \emph{cluster-tilted algebra} if it is the relation extension of a tilted algebra. 

\begin{remark}
 Originally, cluster-tilted algebras were introduced as bound quiver algebras associated to triangulations of polygons in \cite{CCS}, and as endomorphism algebras of cluster-tilting objects in cluster categories in \cite{BMR}.
\end{remark}
\begin{remark}
Relation extensions which are not necessarily cluster-tilted have been studied in \cite{Am,BFPPT}, and a generalization of the construction is given in \cite{AGS}.
\end{remark}
The relation between a tilted algebra and its cluster-tilted algebra naturally provides a strong connection between the corresponding module categories. For example,
a slice in the module category of a tilted algebra embeds as a local slice in the module category of the corresponding cluster-tilted algebra. 
\begin{thm}
 \cite{ABS2} \label{abs2}
 Let $B$ be a cluster-tilted algebra. Then \begin{enumerate}
\item $B$ admits a local slice.
\item For each local slice $\zS$ in $\zG(\textup{mod}\,B)$, the algebra $C= B/\Ann_B\zS$ is tilted with slice $\zS$ and $B= R(C)$.
\item For every tilted algebra $C$ such that $B=R(C)$ there is a local slice $\zS$ such that $C=B/\Ann_B\zS$. 
\end{enumerate}
\end{thm}

\subsection{One-point extensions}
Let $\zL$ be a $k$-algebra and $M\in \textup{mod}\,\zL$. Then the one-point extension  of $A$ by $M$ is the triangular matrix algebra 
\[ 
\zL[M]=\left[
\begin{array}{ccc}
\zL  & 0   \\
M  & k   
\end{array}\right], 
\textup{ whose elements are of the form }
\left[
\begin{array}{ccc}
a  & 0   \\
m  & \mu   
\end{array}\right] ,\]
with $a\in \zL, m\in M$ and $\mu\in k$; and its multiplication  is given by
\[ 
\left[
\begin{array}{ccc}
a  & 0   \\
m  & \mu   
\end{array}\right] 
\left[
\begin{array}{ccc}
a'  & 0   \\
m'  & \mu'   
\end{array}\right] =
\left[
\begin{array}{ccc}
aa'  & 0   \\
m a'+\mu m' & \mu\mu'   
\end{array}\right] \]
using the right $\zL$-module structure of $M$ as well as its $k$-vector space structure.

The quiver $Q_\zL$ is a full subquiver of $Q_{\zL[M]}$, and $Q_{\zL[M]}$ has exactly one vertex more than $Q_\zL$, the \emph{extension vertex}. Moreover the extension vertex is a source in $Q_{\zL[M]}$ and the radical of the indecomposable projective $\zL[M]$-module at the extension  vertex is equal to $M$.  
The arrows in $Q_{\zL[M]}$ starting at the extension vertex are called \emph{extension arrows}.

In this article, we shall be interested in the special case where the module $M$ is a projective $\zL$-module. In this case, the relation ideal of the one-point extension $\zL[M]$ can be identified with the relation ideal of $\zL$, more precisely, if $\zL=kQ_\zL/I $ where $I$ is an admissible ideal generated by a set of paths in $Q_\zL$ then $\zL[M]=kQ_{\zL[M]}/ I'$ where $I'$ is the admissible ideal in $kQ_{\zL[M]} $ generated by the same set of paths. 
Moreover, for each indecomposable direct summand $P_C(j)$ of the projective module $M$, there is precisely one extension arrow from the extension vertex to the vertex $j$, counting multiplicities. %
%
%
%
%
%
%
%
\section{Main results}\label{sect3}

\subsection{Relation extensions and one-point extensions}
Throughout this subsection, $C$ is a $k$-algebra of global dimension at most 2, $P$ is a projective $C$-module, and we   consider the one-point extension $C[P]$. We   use the notation $1,2,\ldots ,n$ for the vertices of the quiver of $C$,  and we use $n+1$ to denote the extension vertex of $C[P]$. To simplify the notation, we denote the indecomposable projective (respectively injective) $C[P]$-modules simply by $P(i)$ (respectively $I(i)$). For the indecomposable projective (respectively injective) $C$-modules we continue to write $P_C(i)$ (respectively $I_C(i)$).

Let $a_j$ be the multiplicity of the indecomposable projective $C$-module $P_C(j)$ at the vertex $j$ in the direct sum decomposition of $P$, thus
\[P=\bigoplus_{j\in J} a_j P_C(j).\]
Now consider the relation extension $R(C)$ of $C$ and let $\Pbar$ be the projective $R(C)$- module whose indecomposable summands $P_{R(C)}(j)$ have the same multiplicity $a_j$, thus
\[\Pbar=\bigoplus_{j\in J} a_j P_{R(C)}(j).\]

\begin{thm}\label{thm 1}
Let $C$ be a  $k$-algebra of global dimension at most 2, and let $P$ be a projective $C$-module. Let $\Pbar$ be the corresponding projective $R(C)$-module defined above. Then there exists an isomorphism of $k$-algebras
\[R(C[P])\cong R(C)[\Pbar].\]
\end{thm}

\begin{remark} The idea of the proof is the following. 
 It has been shown in \cite{ABS} that if $R(C)$ is the relation extension of an algebra $C$ of global dimension at most 2, then the set of new paths in the quiver of $R(C)$ from a vertex $x$ to a vertex $y$ is given by a basis of $\Ext^2_C(I_C(y),P_C(x))$. Since the one-point extensions in the theorem are with respect to  projective modules, the extension vertex is not involved in any relations, and, consequently, the one-point extension commutes with the relation extension. 
\end{remark}

Before  proving the Theorem \ref{thm 1} we  need the following lemmas.

The $C[P]$-modules can be considered as $C$-modules by restriction of scalars. On the other hand, although $1_C\ne 1_{C[P]}$,  the $C$-modules can be considered as $C[P]$-modules by extending the scalars in the following way. Viewing the modules as quiver representations, extending the scalars just means that the vector space on the extension vertex is the zero space.  

\begin{lem}
 \label{lem11}
 There is an isomorphism of right $C$-modules
 $\Pbar\cong P\oplus \Ext^2_{C}(DC,P)$.
\end{lem}
\begin{pf}
 In $\textup{mod}\,R(C)$, there is a short exact sequence of the form 
 \[0\longto \Ext^2_C(DC,P)\longto \Pbar\longto P\longto 0.\] This shows that there is a $k$ vector space isomorphism $\Pbar\cong P\oplus \Ext^2_{C}(DC,P)$, which is also a $C$-module isomorphism.
 \qed
\end{pf}

\begin{lem}
 \label{lem12}\begin{enumerate}
\item  There is an isomorphism of $C$-bimodules 
\[\psi : \Ext^2_{C[P]}(DC[P],C ) \longto \dcc \]
acting trivially on the new paths in $R(C[P])$ which do not start at the extension vertex $n+1$.
\item  There is an isomorphism of $C$-bimodules 
\[\chi : \Ext^2_{C[P]}(DC[P],P(n+1) ) \longto \Ext^2_C(DC,P) \]
acting on the new paths in $R(C[P])$ which start  at the extension vertex $n+1$ by deleting the initial arrow. 
\item For all $m\in P_c$ and all $g\in\Ext^2_{C[P]}(DC[P],C ) $, we have \[ m\psi(g)=\chi(mg).\]
\end{enumerate}
\end{lem}

\begin{pf}  First note that $I(n+1)$ is a simple module and that
 \begin{equation}\label{ses}0\longto P \longto P(n+1)\longto I(n+1)\longto 0\end{equation}
 is a projective resolution in $\textup{mod}\,C[P]$, and thus, $I(n+1)$ has projective dimension 1, whence
$\Ext^2_{C[P]}(I(n+1),C)=0$.

Now let $i$ be a vertex different from $n+1$. Then the injective $C$-module
$I_C(i)$ considered as a $C[P]$-module is a submodule of the corresponding injective $C[P]$-module $I(i)$, and the quotient is a direct sum of copies of the simple injective module $I(n+1)$. Thus we have a short exact sequence in $\textup{mod}\,C[P]$ of the form
\[0\longto I_C(i)\longto I(i)\longto m I(n+1)\longto 0\]
for some $m\in \mathbb{Z}.$

Applying $\Hom_{C[P]}(-,C) $ to this sequence and using the fact that $\pd I(n+1)=1$ yields an isomorphism
\begin{equation}\label{eq one}\Ext^2_{C[P]} (I(i),C)\longto \Ext^2_{C[P]}(I_C(i),C) .\end{equation}
Moreover, given a minimal projective resolution of $I_C(i)$ in $\textup{mod}\,C$ one obtains a minimal projective resolution of $I_C(i)$ in $\textup{mod}\,C[P]$ simply by extending scalars, because the extension vertex $n+1$ is a source in the quiver of $C[P]$, which is not in the support of $I_C(i)$. Thus there is an isomorphism of left $C$-modules
\[\Ext^2_{C[P]} (I_C(i),C)\longto \Ext^2_{C}(I_C(i),C) .\] 
Composing this isomorphism with the isomorphisms in (\ref{eq one}), we obtain an isomorphism of left $C$-modules 
\[\Ext^2_{C[P]} (I(i),C)\longto \Ext^2_{C}(I_C(i),C) ,\] 
which sends a new path in $R(C[P])$ with terminal point $i$ and starting point  different from $n+1$ to the same path in $R(C)$.
This isomorphism induces an isomorphism $\psi$ of the left $C$-modules \[\Ext^2_{C[P]} (DC[P],C)\longto \Ext^2_{C}(DC,C) .\] 

For the right $C$-module structure, let
\[0\longto P_C(i)\stackrel{f_0}{\longto} I_0 \stackrel{f_1}{\longto} I_1 \stackrel{f_2}{\longto} I_2 \longto 0\]
be a minimal injective resolution in $\textup{mod}\,C$. The indecomposable components of the maps $f_i$ are given by comultiplication of paths in $Q_C$. Then, since $n+1$ is a source in $Q_{C[P]}$, there is a minimal injective resolution in $\textup{mod}\,C[P]$ of the form 
\[0\longto P_C(i)\stackrel{\bar f_0}{\longto} \bar I_0 \stackrel{\bar f_1}{\longto} \bar I_1\oplus a \bar I(n+1)\stackrel{[\bar f_2\ 0]}{\longto} \bar I_2 \longto 0\]
where each indecomposable injective $C[P]$-module $I(j)$ appears in the direct sum decomposition of $\bar I_h$, $h=0,1,2$, with the same multiplicity as the corresponding indecomposable injective $C$-module $I_C(j)$ in $I_h$, and each indecomposable component of $\bar f_h$ is given by the comultiplication with the same path as the corresponding indecomposable component of the map $f_h$.
It follows that $\psi $ is also an isomorphism of right $C$-modules.
This shows 1.

Since the projective dimension of $ I(n+1)$ is 1, we have $\Ext^2_{C[P]}(I(n+1),P(n+1))=0$. Moreover, by a similar argument as in the proof of 1, we get an isomorphism 
\[\Ext^2_{C[P]} (I(i),P(n+1))\longto \Ext^2_{C[P]}(I_C(i),P(n+1)),\]
for all $i\ne n+1$.

Applying the functor $\Hom_{C[P]}(I_C(i),-)$ to the short exact sequence (\ref{ses})
and using the fact that $I(n+1)$ is injective yields
an isomorphism
\[\Ext^2_{C[P]} (I_C(i),P)\longto \Ext^2_{C[P]}(I_C(i),P(n+1)) ,\]
which sends a new path $w$ of $R(C[P])$ starting at a vertex $j$ which corresponds to an indecomposable summand $P_C(j)$ of $P$, to the new path $\alpha w$ in $R(C[P])$, where $\alpha$ is the extension arrow $(n+1) \to j$ in $C[P]$ corresponding to the indecomposable summand $P_C(j)$ of $P$.

By the same argument as in the proof of 1, there is an isomorphism of vector spaces
\[\Ext^2_{C[P]} (I_C(i),P)\longto \Ext^2_{C}(I_C(i),P) ,\] 
and composing the three isomorphisms yields an isomorphism of vector spaces 
\[\chi : \Ext^2_{C[P]}(DC[P],P(n+1) ) \longto \Ext^2_C(DC,P) \] 
which acts on the new paths in $R(C[P])$ that start  at the extension vertex $n+1$ by deleting the initial arrow.
Finally, by the same argument as in the proof of 1, we see that $\chi $ is an isomorphism of right $C$-modules. 
This show 2. 

It remains to show 3.
For the action of $m$ on $g$, we must consider $m\in C[P]$ as an element of $P=\rad P(n+1)$, thus $m$ is a linear combination of non-constant paths starting at $n+1$. Then $mg$ is given by composing the paths in $m$ and $g$. Applying $\chi$  to $mg$ means deleting the initial arrow from each of these paths, so that the new starting point is in the top of $P$.

On the other hand, for the action of $m$ on $\psi(g)$, we must consider $m\in C$ as an element of the the projective $C$-module $P$, thus $m$ is given by a linear combination of paths in $Q_C$, which are obtained by deleting the initial arrow from each of the paths in the expression for $m \in \rad P(n+1)$ discussed above. 
The product $m\psi(g)$ is then given by the composition of the paths in $m$ and $\psi(g)$. This shows 3. 
\qed
\end{pf}
\smallskip

Proof of Theorem \ref{thm 1}. With the notation of Lemmas \ref{lem11} and \ref{lem12}, we define a map
\[\phi\colon R(C[P]) \longto R(C)[\Pbar_{R(C)}]\] by 
\[ \phi\left(\left[
\begin{array}{ccc}
 c & 0     \\
 m &   \mu 
\end{array}   
\right] ,
\left[
\begin{array}{ccc}
 g & 0     \\
 h &  0
\end{array} \right]
\right)
\quad =\quad
\left[
\begin{array}{ccc}
 (c,\psi(g)) & 0     \\
( m , \chi(h) )&   \mu 
\end{array} \right]
\]
where $c\in C,\, m\in P, \,\mu \in k$, so $\left[
\begin{array}{ccc}
 c & 0     \\
 m &   \mu 
\end{array}   
\right]  $
is an element of $C[P]$, and $g\in \Ext^2_{C[P]}(DC[P], C)$, $h \in \Ext^2_{C[P]}(DC[P], P(n+1))$, and the two zeros in the matrix $\left[
\begin{array}{ccc}
 g & 0     \\
 h &  0
\end{array} \right]$
correspond to  $\Ext^2_{C[P]}(I(n+1) ,C[P])=0$. The map $\phi $ is clearly bijective and the following computation shows that $\phi$ is a homomorphism of $k$-algebras.

First note that 
\[\left(
\left[
\begin{array}{ccc}
 c & 0     \\
 m &   \mu 
\end{array}   
\right] ,
\left[
\begin{array}{ccc}
 g & 0     \\
 h &  0
\end{array} \right]
\right)
\left(\left[
\begin{array}{ccc}
 c' & 0     \\
 m' &   \mu' 
\end{array}   
\right] ,
\left[
\begin{array}{ccc}
 g' & 0     \\
 h' &  0
\end{array} \right]
\right)
\] is equal to 
\[\left(\left[
\begin{array}{ccc}
 cc' & 0     \\
 mc'+\mu m' &   \mu \mu'
\end{array}   
\right] ,
\left[
\begin{array}{ccc}
c g'+gc' & 0     \\
mg'+\mu h' +hc' &  0
\end{array} \right]\right).
\]
Applying $\phi$ to this expression yields
\[\left[
\begin{array}{ccc}
 (cc',\psi(cg'+gc') )& 0     \\
 (mc'+\mu m', \chi( mg'+\mu h' +hc' ) )&   \mu \mu'
\end{array}   
\right] .\]

On the other hand, 
\[\phi\left(
\left[
\begin{array}{ccc}
 c & 0     \\
 m &   \mu 
\end{array}   
\right] ,
\left[
\begin{array}{ccc}
 g & 0     \\
 h &  0
\end{array} \right]
\right)\phi
\left(\left[
\begin{array}{ccc}
 c' & 0     \\
 m' &   \mu' 
\end{array}   
\right] ,
\left[
\begin{array}{ccc}
 g' & 0     \\
 h' &  0
\end{array} \right]
\right)
\] 
is equal to 
\[\left[
\begin{array}{ccc}
 (cc',c\psi(g')+\psi(g)c' )& 0     \\
 (mc'+\mu m', \mu\chi(  h') +\chi(h)c'  +m\psi(g') )&   \mu \mu'
\end{array}   
\right] \]
and the result follows from Lemma \ref{lem12}.
\qed 

\subsection{One-point extensions of cluster-tilted algebras}
We now study the question when the one-point extension of a cluster-tilted algebra is again cluster-tilted. The following lemma seems to be well-known; we provide a proof for the convenience of the reader.
%
%
%
%
%
 
 \begin{lem}\label{lemneu}
Let $C$ be a tilted algebra with slice $\Sigma$ and let $M$ be a module in $\Sigma$. Then the one-point extension $C[M]$ is tilted with slice $\Sigma'=\Sigma \cup P(n+1)$.
\end{lem}
\begin{pf}
 Let $\{1,2, \ldots,n\}$ be the set of vertices in $Q_C$ and let $n+1$ denote the extension vertex in $Q_{C[M]}$. The full subquivers in both Auslander-Reiten quivers $\zG(\textup{mod}\,C)$ and $\zG(\textup{mod}\,C[M])$, whose  points are the predecessors of $\zS$, are equal. In $\zG(\textup{mod}\,C[M])$, the new indecomposable projective $C[M]$-module $P(n+1)$ lies on a new $\tau_{C[M]}$-orbit with one arrow $M'\to P(n+1)  $ for each indecomposable summand $M'$ of $M$, counting multiplicities.
Therefore, the full subquiver given by $\zS'=\zS\cup\{P(n+1)\}$ satisfies the condition (s1) and (s2') in the definition of a left section.

To show that $\zS'$ satisfies condition (s3), observe that any non-constant path in $\zG(\textup{mod}\,C[P])$ that ends in $P(n+1)$ must pass through one of its immediate predecessors, which are given by the indecomposable summands of $M$, because $M$ is the radical of $P(n+1)$. Since  $M$ lies in $\zS$, it follows that every path in $\zG(\textup{mod}\,C[M])$ that ends in $P(n+1)$ restricts to a path that ends in a point in $\zS$. On the other hand, every path in $\zG(\textup{mod}\,C[M])$ that ends in a point in $\zS$ is actually a path in $\zG(\textup{mod}\,C)$. Hence the condition (s3) for $\zS'$ in  $\zG(\textup{mod}\,C[M])$ follows from the condition  (s3) for $\zS$ in  $\zG(\textup{mod}\,C)$. This shows that $\zS' $ is a left section.

Our next goal is to show that $\Hom_{C[M]}(\tau_{C[M]}^{-1}E',E'')=0$ for all $E',E''\in \zS'$. Clearly it suffices to check this property for indecomposable modules $E',E''$. Let us suppose first that $\Hom_C(M,E')=0$. In that case it follows from \cite[Corollary XV 1.7]{SiSk} that $\tau^{-1}_{C[M]} E'=\tau^{-1}_{C} E'$. Thus if $E''\in \zS$ then $\Hom_{C[M]}(\tau_{C[M]}^{-1}E',E'')=0$ because $\zS$ is a slice in $\textup{mod}\,C$. On the other hand, if $E''=P(n+1)$ then any non-zero morphism in $\Hom_{C[M]}(\tau_{C[M]}^{-1}E',E'')$ would factor through $M=\rad P(n+1) \in \zS$ and, again, we conclude that $\Hom_{C[M]}(\tau_{C[M]}^{-1}E',E'')=0$.

Now suppose that $\Hom_C(M,E')\ne 0$. It follows from the construction of $\zS'$  that $E'$ is either a direct summand of $M$ or $E'=P(n+1)$. Suppose first that $E'$ is a summand of $M$. Without loss of generality, we may assume that $M$ has no injective summands.
 Consider the following almost split sequence in $\textup{mod}\,C[M]$ 
\begin{equation}\label{eqCM}0\longto M\longto P(n+1) \oplus X\longto \tau^{-1}_{C[M]}M\longto 0 .\end{equation}
Then  $X$ is a $C$-module, and   we have the following almost split sequence in $\textup{mod}\,C$
\begin{equation}\label{eqC} 0\longto M\longto   X\longto \tau^{-1}_{C}M\longto 0 .\end{equation}

If  $E''\in \zS$, then applying $\Hom_C(-, E'')$ to the exact sequence (\ref{eqC}) shows that $\Hom_C(X,E'')\cong \Hom _C(M,E'')$. Moreover,  $\Hom( P(n+1) ,E'')=0$, since $E''$ is not supported on the extension vertex $n+1$. Therefore, applying $\Hom(-,E'')$ to the exact sequence (\ref{eqCM}) yields an exact sequence
\[ 0\to  \Hom(\tau^{-1}_{C[M]}M, E'') \stackrel{}{\longto} \Hom(X,E'') \stackrel{g}{\longto} \Hom(M,E''),\] with $g$ an isomorphism. 
 Consequently, 
$ \Hom(\tau^{-1}_{C[M]}M, E'') =0 $ if $E''\in\zS$. 

On the other hand, if $E''=P(n+1)$ then $\Hom(X,E'')=0$. Therefore applying $\Hom(-,E'')$ to the exact sequence (\ref{eqCM}) yields an exact sequence
\[ 0\to  \Hom(\tau^{-1}_{C[M]}M, E'') \stackrel{}{\longto} \Hom(P(n+1),E'') \stackrel{g'}{\longto} \Hom(M,E''),\]
where 
the morphism $g'$ is injective, since $M$ is the radical of $E''=P(n+1)$, and, again, we conclude that $ \Hom(\tau^{-1}_{C[M]}M, E'') =0 $.

It remains the case where $E'=P(n+1)$. 
Since $M$ lies in the slice $\zS$ in $\textup{mod}\,C$, we have $\id\!_CM\le 1$.
Let 
\[0\longto M\longto I^0_C(M)\longto I^1_C(M)\longto0\] 
be a minimal injective resolution of $M$ in $\textup{mod}\,C.$ Then we have a minimal injective resolution in $\textup{mod}\,C[M]$ of the form
\[0\longto M\longto I^0(M)\longto I^1(M)\oplus m I(n+1) \longto0\] 
where $m\ge 0$ and the multiplicity of each indecomposable injective $C[M]$-module $I(j)$ at $ j\ne n+1$ in $I^0(M)$ (respectively $I^1(M)$) is equal to the multiplicity of the indecomposable injective $C$-module $I_C(j)$ in $I_C^0(M)$ (respectively $I_C^1(M)$). In particular, $\id\!_{C[M]} M\le 1$. Using the horseshoe lemma on the short exact sequence 
\[0\longto M\longto P(n+1)\longto I(n+1)\longto 0\]
we conclude that $\id\!_{C[M]}P(n+1)\le 1$, and  a standard Lemma from representation theory, see for example \cite[Lemma IV 2.7]{ASS}, implies \[\Hom(\tau^{-1}_{C[M]} P(n+1), C[M])=0.\] Now
 the Auslander-Reiten formula yields an isomorphism
\[ \Hom(\tau^{-1}_{C[M]} P(n+1),E'')\cong D\Ext^1(E'',P(n+1)).\]
 If $E''=P(n+1)$ then the right hand side is zero. To compute the right hand side in the case where $E''\in\zS$, we  apply the functor $\Hom_{C[M]}(E'',-)$ to the short exact sequence 
\[ 0\longto M\longto  P(n+1) \longto I(n+1)\to 0\]
to conclude that there exists  a surjective map \[\Ext^1_{C[M]}(E'',M)\to \Ext^1_{C[M]}(E'',P(n+1)).\]  
But    $\Ext^1_{C[M]}(E'',M) = 0$, because both $E''$ and $M$ lie in the slice $\zS$ of $\textup{mod}\,C$. 
This completes the proof that  $\Ext^1_{C[M]}(E'',E')=0$, for all $E',E''\in \zS'$.
\smallskip

We have shown that $\zS'$ satisfies the hypotheses of Theorem \ref{thm assem}. Moreover, $\Ann_C\zS=0$, because $\zS$ is a slice in $\textup{mod}\,C$, hence a faithful $C$-module. This implies that  $\Ann_{C[M]}\zS'=0$, because the only paths in the quiver $Q_{C[M]}$ of $C[M]$, which are not already in the subquiver $Q_C$, must start at the extension vertex $n+1$ and therefore do not annihilate the direct summand $P(n+1)$ of $\zS'$. 
Now Theorem \ref{thm assem} implies that $C[M]$ is a tilted algebra with slice $\zS'$.\qed
\end{pf}
\bigskip
 
\begin{thm}\label{thm main}
 Let $B$ be a cluster-tilted algebra and $P$ a projective $B$-module whose indecomposable summands lie in a local slice of $\zG(\textup{mod}\,B)$. Then $B[P]$ is cluster-tilted. 
\end{thm}

\begin{remark}
 Note that the extension vertex in $B[P]$ is a source.
\end{remark}
 
\begin{pf}
Let $\zS$ be a local slice of $B$ containing $P$. By Theorem \ref{abs2}, the algebra $C= B/\Ann \zS$ is   tilted   with slice $\zS$, and $B=R(C)$ is the relation extension of $C$. Now it follows from Lemma \ref{lemneu} that $C[P]$ is a tilted algebra with slice $\zS'=\zS\cup\{P(n+1)\}$,
 and  Theorem \ref{thm 1} implies that the one-point extension $B[P]$ is isomorphic to the relation extension $R(C[P])$ of the tilted algebra $C[P]$. By definition, this means that $B[P]$ is a cluster-tilted algebra.
\qed
\end{pf}

\begin{remark}
 Dually, if $I$ is an injective $B$-module whose indecomposable summands lie in a local slice of $\zG(\textup{mod}\,B)$, then the one-point coextension of $B$ with respect to $I$ is cluster-tilted. In this situation, the extension vertex is a sink.
\end{remark}

\subsection{Further extensions of cluster-tilted algebras}
We consider now a different construction of an extension of a cluster-tilted algebra, which produces a cluster-tilted algebra whose extension vertex can have incoming \emph{and} outgoing arrows.
\begin{thm}\label{thm main2}
 Let $B$ be a cluster-tilted algebra with local slice $\zS$, and let $C=B/\Ann \zS$ be the corresponding tilted algebra. Let $M$ be any module in $\zS$. Then $B'=R(C[M])$ is a cluster-tilted algebra such that
 \begin{itemize}
\item[ \textup{(a)}]
$Q_{B'}$ contains $Q_{B}$ as a full subquiver and $Q_{B'}$ has one new vertex $n+1$;
\item[ \textup{(b)}] $\zS'=\zS\cup\{P(n+1)\}$ is a local slice;
\item[ \textup{(c)}]$B'/B'e_{n+1} B' = B$;
 \item[ \textup{(d)}]
$\rad P(n+1)=M$ and $ I(n+1)/S(n+1)= \tau_B M=\tau_{B'}M$.
\end{itemize}
Moreover, the quiver of $B'$ is obtained from the quiver of $B$ by adding one vertex $x$ and one arrow   $n+1\to i $ for each indecomposable summand $S(i)$ of $\textup{top}\,M$, and one arrow $j\to n+1$ for each indecomposable summand $S(j)$ of $I(n+1)/S(n+1)$, counting multiplicities.
\end{thm}

\begin{pf}
Lemma \ref{lemneu} implies that $C[M]$ is a tilted algebra with slice $\zS'=\zS\cup\{P(n+1)\}$ whose quiver has one more vertex than the quiver of $C$, 
 and therefore $B'$ is a cluster-tilted algebra with local slice $\zS'$ and one new vertex $n+1$, which proves (a) and (b). 

(c) The algebra $B' / B' e_{n+1} B'$ is  cluster-tilted by   \cite[Theorem 2.13]{BMR4}, and since cluster-tilted algebras are determined by their quivers, it suffices to show that the quiver   of $B$ is the same as the quiver  of $B' / B' e_{n+1} B'$. 
The quivers of the relation extensions $B=R(C)$ and $B'=R(C[M])$ are obtained from the quivers of $C$ and $C[M]$ respectively by adding $e_{ij}$ arrows $i\to j$ where $e_{ij}$ is the $k$-dimension of the vector space $\Ext^2_C(S(i), S(j))$ and $\Ext^2_{C[M]}(S(i), S(j))$, respectively. Now if $i $ and $j$ are different from $n+1$ then the dimensions of $\Ext^2_C(S(i), S(j))$ and $\Ext^2_{C[M]}(S(i), S(j))$ are equal. Moreover the quivers of $C$ and $C[M]$ differ only by arrows starting at the extension vertex $n+1$. Since the point $n+1$ is deleted when taking the quotient $B' / B' e_{n+1} B'$, we conclude that the quivers of $B$ and $B' / B' e_{n+1} B'$ are equal. 

(d) Since $P(n+1)$ lies in the local slice $\zS'$ it follows that its radical in $\textup{mod}\,B'$ is the same as its radical in $\textup{mod}\,C[M]$, thus $\rad P(n+1) = M$. 

On the other hand, since the vertex $n+1$ is a source in the quiver of $C[M]$, the socle factor $I(n+1)/S(n+1)$ of the corresponding injective in the relation extension $R(C[M])$ is given by $D\Ext^2_{C[M]}(I(n+1),C[M])$.
 In order to compute the latter, we 
apply $\Hom(-,C[M])$ to the short exact sequence 
\[0\longto M\longto P(n+1)\longto I(n+1)\longto 0\]
to get an isomorphism
\[ \Ext^1_{C[M]}(M,C[M])\cong \Ext^2_{C[M]}(I(n+1),C[M]).\]
 Moreover, since the projective dimension of the $C[M]$-module $M$ is at most 1, the Auslander-Reiten formula implies that 
 \[ D\Hom_{C[M]}(C[M],\tau_{C[M]} M)\cong \Ext^1_{C[M]}(M,C[M]),\] and finally there is an isomorphism $\Hom_{C[M]}(C[M],\tau_{C[M]} M)\cong \tau_{C[M]} M$. 
Combining these three isomorphisms shows that $I(n+1)/S(n+1)=\tau_{C[M]} M$, and (d) follows because $\tau_{C[M]} M =\tau_{C} M =\tau_{B} M =\tau_{B'} M$, since $M$ lies in $\zS$, see \cite[Proposition 3]{ABS4}.

The statement about the arrows in the quiver of $B'$ follows from (d).
\qed
\end{pf}

\begin{remark}
 If $M$ is projective then the cluster-tilted algebra $B'$ of the theorem is the one-point extension $B[P]$ of $B$ in $P$. 
\end{remark}
\begin{remark}
 In any cluster-tilted algebra, $\tau \,\rad P(i)$ is always equal to $I(i)/S(i)$, see for example \cite[Lemma 5]{ABS4}.
\end{remark}

\subsection{Examples} 
\begin{example}
 Let $B$ be the cluster-tilted algebra of type $\mathbb{D}_4$ given by the quiver 
\[\xymatrix{&2\ar[rd]^\zb\\ 
1\ar[ru]^\za\ar[rd]_\zg &&4\ar[ll]_\ze\\
&3\ar[ru]_\zd }\]
bound by the relations $\za\zb+\zg\zd,\,\ze\za, \, \ze\zg,\,\zb\ze,\,\zd\ze$. It's Auslander-Reiten quiver is the following
$$\xymatrix@C=10pt@R=0pt
{
&&{\begin{array}{c} \mathbf 2\\\mathbf 4 \end{array}}\ar[dr]&&
{\begin{array}{c} 3 \end{array}}\ar[dr]&&
{\begin{array}{c} 1\\2 \end{array}}\ar[dr]&&
&&
\cdots
\\
\cdots\ar[r]&{\begin{array}{c} 4 \end{array}}\ar[ur]\ar[dr]&
&
{\begin{array}{c} \mathbf{ 2\ 3}\\\mathbf 4 \end{array}}\ar[dr]\ar[r]\ar[ur]&
{\begin{array}{c} \mathbf 1\\\mathbf{ 2\ 3} \\\mathbf 4\end{array}}\ar[r]&
{\begin{array}{c} 1\\2 \ 3 \end{array}}\ar[dr]\ar[ur]&
&
{\begin{array}{c} 1 \end{array}}\ar[r]&
{\begin{array}{c} 4\\1 \end{array}}\ar[r]&
{\begin{array}{c} 4 \end{array}}\ar[ur]\ar[dr]&
\cdots
\\
&&{\begin{array}{c}\mathbf 3\\\mathbf 4\\ \  \end{array}}\ar[ur]&&
{\begin{array}{c}\ \\ 2\\ \  \end{array}}\ar[ur]&&
{\begin{array}{c}1\\ 3\\ \  \end{array}}\ar[ur]&&
&&
\cdots
}$$
where the two points with label $4$ must be identified. 

To illustrate Theorem \ref{thm main}
let $P=P(1)\oplus P(2) \oplus P(3) $. Then $P$ lies on the local slice $\zS$ which is illustrated by the boldface modules in the Auslander-Reiten quiver above, and the one-point extension $B[P]$ is the cluster-tilted given by the following quiver  bound by the same relations.

\[
\xymatrix{&&&2\ar[rd]^\zb\\ 
5\ar[rr]\ar@/^10pt/[rrru]\ar@/_10pt/[rrrd]&&1\ar[ru]^\za\ar[rd]_\zg &&4\ar[ll]_\ze\\
&&&3\ar[ru]_\zd 
}\]

\bigskip
To illustrate Theorem \ref{thm main2}, let $M=2$ be the simple module $S(2)$. Then $\tau M = P(3)$. The module $M$ lies on a local slice which is obtained from the local slice $\Sigma $ above by replacing $P(3)$ with $M$. The corresponding tilted algebra $C$ is given by the quiver
\[\xymatrix@R10pt{&2\ar[rd]^\zb\\ 
1\ar[ru]^\za\ar[rd]_\zg &&4\\
&3\ar[ru]_\zd }\] 
bound by the relation $\za\zb+\zg\zd$. The one-point extension $C[M]$ is given by the quiver 
\[\xymatrix@R10pt{ &&5\ar[dl]_\rho\\&2\ar[rd]^\zb\\ 
1\ar[ru]^\za\ar[rd]_\zg &&4\\
&3\ar[ru]_\zd }\] 
bound by the relations $\za\zb+\zg\zd, \,\rho\zb.$ Finally the relation extension $B'=R(C[M])$ is  given by the quiver 
\[\xymatrix{ &&5\ar[dl]_\rho\\&2\ar[rd]^\zb\\ 
1\ar[ru]^\za\ar[rd]_\zg &&4\ar[uu]_\zs\ar[ll]_\ze\\
&3\ar[ru]_\zd }\] 
bound by the relations $\beta\epsilon,\, 
\epsilon\alpha+\sigma\rho,\,
\delta\epsilon,\,
\epsilon\gamma,\,
\beta\sigma,\,
\alpha\beta+\gamma\delta
,\,\rho\beta$.
\end{example}

\begin{example}\label{exnonex}
This example shows that the condition that the projective module $P$ lies on a local slice in Theorem \ref{thm main} is necessary, even if $P$ is indecomposable. 
Let $B$ be the cluster-tilted algebra of type $\tilde{\mathbb{A}}_{3,1}$ given by the following quiver bound by the relations $\za\zb,\zb\zg,\zg\za$.
\[\xymatrix{1\ar@<2pt>[dd]^\zb\ar@<-2pt>[dd]
\\
&3\ar[lu]_\za\ar[r]&4\\
2\ar[ru]_\zg
}
\]
To see that $B$ is cluster-tilted, one can apply the mutations in vertex 3 and then in vertex 4 to obtain an acyclic $\tilde{\mathbb{A}}_{3,1}$ quiver.

The projective $B$-module $P(3)$ lies in a regular component of the Auslander-Reiten quiver, and therefore there is no local slice containing $P(3)$. The one-point extension $B[P(3)]$ has the following quiver. \[\xymatrix{1\ar@<2pt>[dd]^\zb\ar@<-2pt>[dd]
\\
&3\ar[lu]_\za\ar[r]&4\\
2\ar[ru]_\zg && 5\ar[lu]
}
\]
This algebra is not cluster-tilted. In fact, this algebra corresponds to the triangulation of a punctured annulus shown in Figure \ref{fig}, which means that its quiver is not mutation equivalent to an acyclic quiver, see \cite{FST}, and therefore the algebra is not cluster-tilted.
\begin{figure}
\begin{center}
\scalebox{0.5}{\includegraphics{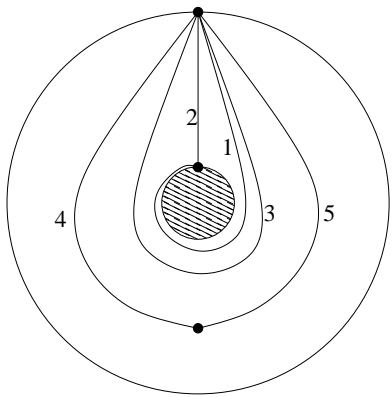}}
\caption{A triangulation of an annulus with one puncture which corresponds to the one-point extension in Example \ref{exnonex}.}\label{fig}
\end{center}
\end{figure}
\end{example}

\end{document}